\documentclass[a4paper,11pt]{amsart}
\usepackage{graphicx}
 \usepackage{mathptmx}
\usepackage{amsmath}
\usepackage{amssymb}

\usepackage{xcolor}

\newtheorem{theorem}{Theorem}[section]

\newtheorem{remark}[theorem]{Remark}

\newcommand{\et}{e^{\frac{t}{2}}-e^{-\frac{t}{2}}}

\begin{document}

\title{Some identities on  $r$-central factorial numbers and $r$-central Bell polynomials
}
\author{ Kim, Dae San }
\address{Department of Mathematics, Sogang University, Seoul 04107, Republic of Korea}
\email{dskim@sogang.ac.kr}

\author{Dolgy, Dmitry V.}
\address{Kwangwoon Institute for Advanced Studies, Kwangwoon University, Seoul 01897, Republic of Korea}
\email{d\underline{ }dol@mail.ru}

\author{Kim, Dojin}
\address{Department of Mathematics, Pusan National University, Busan 46241, Republic of Korea}
 \email{kimdojin@pusan.ac.kr}

 \author{Kim, Taekyun}
\address{Department of Mathematics, Kwangwoon University, Seoul 01897, Republic of Korea}
\email{tkkim@kw.ac.kr}

\maketitle
 
\begin{abstract}
In this paper, we introduce the extended $r$-central factorial numbers of the second and first kinds and the extended $r$-central Bell polynomials, as extended versions and central analogues of some previously introduced numbers and polynomials. Then we study various properties and identities related to these numbers and polynomials and also their connections.
\end{abstract}


\section{Introduction}\label{sec1}
For $n\in \mathbb{N}\cup\{ 0\}$, as is well known, the central factorials $x^{[n]}$ are defined by
\begin{equation}\label{eq1}
  x^{[0]}=1,\quad x^{[n]}=x(x+\frac{n}{2}-1)\cdots (x-\frac{n}{2}+1),\quad (n\geq 1),\quad (\text{see \cite{ref4,ref5,ref6,ref7,ref8,ref15}}).
\end{equation}
It is also well known that the central factorial numbers of the second kind $T(n,k)$ are defined by
\begin{equation}\label{eq2}
  x^{n}=\sum_{k=0}^nT(n,k) x^{[k]}, \quad (\text{see \cite{ref4,ref5,ref15}}),
\end{equation}
where $n$ is a nonnegative integer.

From \eqref{eq2}, we can derive the generating function for $T(n,k),~(0 \leq k\leq n)$ as follows:
\begin{equation}\label{eq3}
 \frac{1}{k!}\left(\et\right)^k=\sum_{n=k}^\infty T(n,k) \frac{t^n}{n!}, \quad (\text{see \cite{ref4,ref5,ref6}}).
\end{equation}

Recently, Kim-Kim \cite{ref8} considered the central Bell polynomials given by
\begin{equation}\label{eq4}
e^{x\left(\et\right)}=\sum_{n=0}^\infty B_n^{(c)}(x)  \frac{t^n}{n!}.
\end{equation}
When $x=1$, $B_n^{(c)}=B_n^{(c)}(1)$ are called the central Bell numbers.

From \eqref{eq4}, we can find the Dobinski-like formula for $B_n^{(c)}(x)$:
\begin{equation}\label{eq5}
 B_n^{(c)}(x) = \sum_{l=0}^\infty \sum_{j=0}^\infty \binom{l+j}{j} (-1)^j \frac{1}{(l+j)!}\left(\frac{l}{2}-\frac{j}{2}\right)^nx^{l+j },\quad  (\text{see \cite{ref8}}).
\end{equation}
The Stirling numbers of the second kind are defined by
\begin{equation}\label{eq6}
 x^n = \sum_{k=0}^n S_2(n,k) (x)_k,~(n\geq 0),\quad  (\text{see \cite{ref1,ref2,ref3,ref9,ref10,ref11,ref16}}),
\end{equation}
where $(x)_0=1,~ (x)_n=x(x-1)(x-2)\cdots(x-n+1),~(n\geq1)$.

From \eqref{eq6}, we can easily derive the following equation \eqref{eq7}.
\begin{equation}\label{eq7}
\frac{1}{k!}(e^t-1)^k= \sum_{n=k}^\infty S_2(n,k)\frac{t^n}{n!},~(k\geq 0),\quad  (\text{see \cite{ref3,ref12,ref13}}).
\end{equation}

In this paper, we introduce the extended $r$-central factorial numbers of the second and first kinds and the extended $r$-central Bell polynomials, and study various properties and identities related to these numbers and polynomials and their connections.
The extended $r$-central factorial numbers of the second kind are an extended version of the central factorial numbers of the second kind and also a ‘central analogue’ of the $r$-Stirling numbers of the second; the extended $r$-central Bell polynomials are an extended version of the central Bell polynomials and also a central analogue of $r$-Bell polynomials; the extended $r$-central factorial numbers of the first kind are an extended version of the central factorial numbers of the first kind and a central analogue of the (unsigned) $r$-Stirling numbers of the first kind. All of these numbers and polynomials were studied before (see \cite{ref4,ref5,ref8,ref14,ref14-1,ref17}).

\section{Extended $r$-central factorial numbers of the second kind and extended $r$-central Bell polynomials}\label{sec2}

Let us first note that, by \eqref{eq3} and \eqref{eq4},
\begin{equation}\label{eq8}
\begin{aligned}
\sum_{n=0}^\infty B_n^{(c)}(x)  \frac{t^n}{n!}&= \sum_{k=0}^\infty x^k\frac{1}{k!}\left(\et\right)^k\\
&=\sum_{k=0}^\infty x^k \sum_{n=k}^\infty T(n,k)\frac{t^n}{n!}\\
&=\sum_{n=0}^\infty\left( \sum_{k=0}^n x^k T(n,k)\right)\frac{t^n}{n!}.
\end{aligned}
\end{equation}
Comparing the coefficients on both sides of \eqref{eq8}, we have
\begin{equation}\label{eq9}
  B_n^{(c)}(x) =  \sum_{k=0}^n x^k T(n,k),\quad (n\geq0).
\end{equation}
For any nonnegative integer $r$, we introduce the extended $r$-central factorial numbers of the second kind given by
\begin{equation}\label{eq10}
\frac{1}{k!}e^{rt} \left( \et\right)^k=\sum_{n=k}^\infty T_{r}(n+r,k+r)\frac{t^n}{n!},\quad (k\geq0).
\end{equation}

\begin{remark}
In \cite{ref6}, the extended central factorial numbers of the second kind were defined as
\begin{equation*}
\frac{1}{k!}(e^{\frac{t}{2}} - e^{-\frac{t}{2}} +rt)^k = \sum_{n=k}^\infty T^{(r)}(n,k) \frac{t^n}{n!}.
\end{equation*}
Note that these numbers are different from the extended $r$-central factorial numbers of the second kind defined in \eqref{eq10}.
\end{remark}

From \eqref{eq3} and \eqref{eq10}, we see that
\begin{equation}\label{eq11}
\begin{aligned}
\sum_{n=k}^\infty T_r(n+r,k+r)\frac{t^n}{n!}&= \frac{1}{k!}e^{rt}\left(\et\right)^k\\
&=\left(\sum_{l=k}^\infty T(l,k)\frac{t^l}{l!}\right)\left(\sum_{m=0}^\infty r^m\frac{t^m}{m!}\right)\\
&=\sum_{n=k}^\infty\left( \sum_{l=k}^n \binom{n}{l} T(l,k)r^{n-l} \right)\frac{t^n}{n!}.
\end{aligned}
\end{equation}
Therefore, by comparing the coefficients on both sides of \eqref{eq10}, the following identity holds.
\begin{theorem}\label{thm1}
  For $n,k,r\in\mathbb{N}\cup\{0\}$ with $n\geq k$, we have
  $$
  T_r(n+r,k+r)=\sum_{l=k}^n\binom{n}{l}T(l,k)r^{n-l}.
  $$
\end{theorem}

Next, we write $e^{(r+x)t}$ as follows:
\begin{equation}\label{eq12}
\begin{aligned}
e^{(r+x)t}&=e^{rt}\sum_{k=0}^\infty (x)_k \frac{1}{k!}\left(\et\right)^k e^{\frac{1}{2}kt}\\
&=\sum_{k=0}^\infty (x)_k \sum_{l=k}^\infty T_r(l+r,k+r)\frac{t^l}{l!}\sum_{j=0}^\infty \left(\frac{k}{2}\right)^j\frac{t^j}{j!}\\
&=\sum_{k=0}^\infty(x)_k \sum_{n=k}^\infty\left( \sum_{l=k}^n \binom{n}{l} T_r(l+r,k+r)\left(\frac{k}{2}\right)^{n-l} \right)\frac{t^n}{n!}\\
&=\sum_{n=0}^\infty \left(\sum_{k=0}^n \sum_{l=k}^n(x)_k \binom{n}{l} T_r(l+r,k+r)\left(\frac{k}{2}\right)^{n-l} \right)\frac{t^n}{n!}.
\end{aligned}
\end{equation}
On the other hand, $e^{(r+x)t}$ can be written as
\begin{equation}\label{eq13}
 e^{(r+x)t}= \sum_{n=0}^\infty (r+x)^n \frac{t^n}{n!}.
\end{equation}

Therefore, by the two expressions in \eqref{eq12} and \eqref{eq13} for  $e^{(r+x)t}$, we obtain the following identity.

\begin{theorem}\label{thm2}
For $n\geq 0$, we have
$$
(r+x)^n=\sum_{l=0}^n \sum_{k=0}^l(x)_k \binom{n}{l} T_r(l+r,k+r)\left(\frac{k}{2}\right)^{n-l}.
$$
\end{theorem}

In view of \eqref{eq4}, we may now introduce the extended $r$-central Bell polynomials associated with the extended $r$-central factorial numbers of the second kind given by
\begin{equation}\label{eq14}
 e^{rt} e^{x\left(\et\right)}= \sum_{n=0}^\infty B^{(c,r)}_n(x) \frac{t^n}{n!}.
\end{equation}

\begin{remark}
In \cite{ref6}, the extended central Bell polynomials were defined as
\begin{equation*}
e^{x (e^{\frac{t}{2}} - e^{-\frac{t}{2}}+rt)} = \sum_{n=0}^\infty Bel_n^{(c,r)}(x) \frac{t^n}{n!}, \,\,\, \left(r \in \mathbb{R}\right).
\end{equation*}
Observe here that these polynomials are different from the extended $r$-central Bell polynomials in \eqref{eq14}.
\end{remark}

From \eqref{eq14}, we note that
\begin{equation}\label{eq15}
\begin{aligned}
\sum_{n=0}^\infty B^{(c,r)}_n(x) \frac{t^n}{n!}&=e^{rt}e^{x \left(\et\right)}\\
&=\sum_{k=0}^\infty x^k \frac{1}{k!}e^{rt}\left(e^{\frac{t}{2}}-e^{-\frac{t}{2}}\right)^k\\
&= \sum_{k=0}^\infty x^k \sum_{n=k}^\infty T_r(n+r,k+r)\frac{t^n}{n!}\\
&=\sum_{n=0}^\infty \left( \sum_{k=0}^n x^k T_r(n+r,k+r)\right)\frac{t^n}{n!}.
\end{aligned}
\end{equation}
By the comparison of the coefficients on both sides of \eqref{eq15}, we can establish the following theorem.
\begin{theorem}\label{thm3}
  For $n\geq 0$, we have that
  $$
  B^{(c,r)}_n(x) =\sum_{k=0}^n x^k T_r(n+r,k+r).
  $$
\end{theorem}

Next, we observe that
\begin{equation}\label{eq16}
\begin{aligned}
\frac{1}{k!}e^{rt}\left(\et\right)^k&=\frac{1}{k!}\sum_{l=0}^k \binom{k}{l}(-1)^{k-l}e^{(l+r-\frac{k}{2})t}\\
&=\sum_{n=0}^\infty \left(\frac{1}{k!}\sum_{l=0}^k \binom{k}{l}(-1)^{k-l} \left(l+r-\frac{k}{2}\right)^n\right) \frac{t^n}{n!}.
\end{aligned}
\end{equation}
By using the central difference operator $\delta$, which is defined by
\begin{equation}\label{eq17}
  \delta f(x)=f\left(x+\frac{1}{2}\right)-f\left(x-\frac{1}{2}\right),
\end{equation}
we can show that
\begin{equation}\label{eq18}
  \delta^k f(x)=\sum_{l=0}^k \binom{k}{l}(-1)^{k-l}f\left(x+l-\frac{k}{2}\right), \quad (k\in \mathbb{N}\cup\{0\}).
\end{equation}
We combine \eqref{eq18} with \eqref{eq16} to derive an equation for $e^{(r+x)t}$  as follows:
\begin{equation}\label{eq19}
\frac{1}{k!}e^{rt}\left(\et\right)^k=\sum_{n=0}^\infty \frac{1}{k!}(\delta^k r^n) \frac{t^n}{n!}.
\end{equation}
From \eqref{eq10} and \eqref{eq19}, we note that
\begin{equation}\label{eq20}
\frac{1}{k!}(\delta^k r^n) =
\begin{cases}
T_r(n+r,k+r), \quad &\text{if}~ n\geq k,\\
0 , &\text{if}~ n<k.
\end{cases}
\end{equation}
Therefore, by \eqref{eq20}, we obtain the following theorem.
\begin{theorem}\label{thm4}
For $n,k\geq 0$, we have
\begin{equation*}
\frac{1}{k!}(\delta^k r^n) =
\begin{cases}
T_r(n+r,k+r), \quad &\text{if}~ n\geq k,\\
0 , &\text{if}~ n<k.
\end{cases}
\end{equation*}
\end{theorem}
By combining Theorems \ref{thm3} and \ref{thm4}, we easily get
\begin{equation}\label{eq21}
  B^{(c,r)}_n(x)=\sum_{k=0}^n x^k \frac{1}{k!} (\delta^k r^n).
\end{equation}
From \eqref{eq14}, we have
\begin{equation}\label{eq22}
\begin{aligned}
\sum_{n=0}^\infty B^{(c,r)}_n(x)\frac{t^n}{n!}&= e^{rt}e^{x\left(\et\right)}\\
&=\left(\sum_{l=0}^\infty B^{(c)}_l (x) \frac{t^l}{l!}\right)\left(\sum_{m=0}^\infty r^m\frac{t^m}{m!}\right)\\
&=\sum_{n=0}^\infty\left( \sum_{l=0}^n \binom{n}{l} B^{(c)}_{n-l} (x) r^l\right)\frac{t^n}{n!}
\end{aligned}
\end{equation}
Therefore,  by comparing the coefficients on both sides of \eqref{eq22}, we get the following identity.
\begin{theorem}\label{thm5}
 For $n\geq 0$, we have
 $$
 B^{(c,r)}_n(x)=\sum_{l=0}^n \binom{n}{l} B^{(c)}_{n-l} (x) r^l.
 $$
\end{theorem}
By \eqref{eq14}, it can be checked that
\begin{equation}\label{eq23}
\begin{aligned}
\sum_{n=0}^\infty B^{(c,r)}_n(x)\frac{t^n}{n!}&= e^{rt}e^{x\left(\et\right)}\\
&=\sum_{m=0}^\infty x^m \frac{1}{m!}(e^t-1)^me^{(r-\frac{m}{2})t}\\
&=\sum_{m=0}^\infty x^m\sum_{l=m}^\infty S_2(l,m)\frac{t^l}{l!} \sum_{k=0}^\infty \left(r-\frac{m}{2}\right)^k\frac{t^k}{k!}\\
&=\sum_{m=0}^\infty x^m\sum_{n=m}^\infty \left(\sum_{l=m}^n \binom{n}{l} S_2(l,m) \left(r-\frac{m}{2}\right)^{n-l}\right)\frac{t^n}{n!}\\
&=\sum_{n=0}^\infty\left(\sum_{m=0}^n\sum_{l=m}^n x^m\binom{n}{l} S_2(l,m) \left(r-\frac{m}{2}\right)^{n-l}\right)\frac{t^n}{n!}.
\end{aligned}
\end{equation}
Therefore, by comparing the coefficients on both sides of \eqref{eq23}, we establish the following theorem.
\begin{theorem}\label{thm6}
For $n\geq 0$, we have
$$
B^{(c,r)}_n(x)=\sum_{l=0}^n\sum_{m=0}^l x^m\binom{n}{l} S_2(l,m) \left(r-\frac{m}{2}\right)^{n-l}.
$$
\end{theorem}
Now, we observe that
\begin{equation}\label{eq24}
\begin{aligned}
\frac{1}{m!}e^{rt}\left(\et\right)^m \frac{1}{k!}\left(\et\right)^k
&=\frac{(m+k)!}{m!k!}\frac{1}{(m+k)!}e^{rt}\left(\et\right)^{m+k}\\
&=\binom{m+k}{m} \sum_{n=m+k}^\infty T_r(n+r,m+k+r)\frac{t^n}{n!}.
\end{aligned}
\end{equation}
On the other hand, it can be seen that
\begin{equation}\label{eq25}
\begin{aligned}
\frac{1}{m!}e^{rt}\left(\et\right)^m \frac{1}{k!}\left(\et\right)^k&=\left(\sum_{l=m}^{\infty}T_r(l+r,m+r)\frac{t^l}{l!}\right)
\left(\sum_{j=k}^{\infty}T(j,k)\frac{t^j}{j!}\right)\\
&=\sum_{n=m+k}^{\infty} \left(\sum_{l=m}^{n-k}\binom{n}{l}
T_r(l+r,m+r)T(n-l,k)\right) \frac{t^n}{n!}.
\end{aligned}
\end{equation}
Therefore, by \eqref{eq24} and \eqref{eq25}, we obtain the following theorem.
\begin{theorem}\label{thm7}
  For $m,n,k\geq0$ with $n\geq m+k$, we have
  $$
  \binom{m+k}{m}T_r(n+r,m+k+r)=\sum_{l=m}^{n-k}\binom{n}{l}
T_r(l+r,m+r)T(n-l,k).
  $$
\end{theorem}

It is known that the generating function of central factorial is given by
\begin{equation}\label{eq26}
\left( \frac{t}{2}+\sqrt{\frac{1}{4}t^2+1}\right)^{2x}=\sum_{n=0}^\infty x^{[n]}\frac{t^n}{n!},\quad \text{(see \cite{ref4,ref5})}.
\end{equation}
If we let $f(t)=2 \log \left(\frac{t}{2}+\sqrt{1+\frac{t^2}{4}}\right)$, then we can easily show that
\begin{equation}\label{eq27}
 f^{-1}(t)=\et.
\end{equation}
By the simple computations with the expressions in \eqref{eq1} and \eqref{eq2}, we can check that  $e^{(x+r)t}$ can be expressed as follows:
\begin{equation}\label{eq28}
\begin{aligned}
 e^{(x+r)t}&=e^{rt}e^{xt}\\
 &=e^{rt}e^{2x\log\left( \frac{e^{\frac{t}{2}}-e^{-\frac{t}{2}}}{2}+\sqrt{1+\frac{(e^{\frac{t}{2}}-e^{-\frac{t}{2}})^2}{4}}\right)}\\
 &=e^{rt}e^{\log \left( \frac{e^{\frac{t}{2}}-e^{-\frac{t}{2}}}{2}+\sqrt{1+\frac{(e^{\frac{t}{2}}-e^{-\frac{t}{2}})^2}{4}} \right)^{2x}}\\
 &=e^{rt} \left( \frac{e^{\frac{t}{2}}-e^{-\frac{t}{2}}}{2}+\sqrt{1+\frac{1}{4} (\et)^2} \right)^{2x}\\
 &=e^{rt} \sum_{k=0}^\infty x^{[k]} \frac{1}{k!}\left(\et\right)^k\\
 &=\sum_{k=0}^\infty x^{[k]} \sum_{n=k}^\infty T(n+r,k+r)\frac{t^n}{n!} \\
 &=\sum_{n=0}^\infty \left( \sum_{k=0}^n T(n+r,k+r)x^{[k]}\right)\frac{t^n}{n!}.
 \end{aligned}
\end{equation}

Alternatively, the term $e^{(x+r)t}$ is also represented by
\begin{equation}\label{eq29}
e^{(x+r)t}=\sum_{n=0}^\infty (x+r)^n\frac{t^n}{n!}.
\end{equation}
Therefore, by \eqref{eq28} and \eqref{eq29}, the following identity is obtained.

\begin{theorem}\label{thm8}
For $n\geq 0$, we have the following identity
  \begin{equation}\label{eq30}
    (x+r)^n=\sum_{k=0}^n T_r(n+r,k+r) x^{[k]}.
  \end{equation}
\end{theorem}

\section{Extended $r$-central factorial numbers of the first kind }\label{sec3}
Throughout this section, we assume that $r$ is any real number.
The (unsigned) $r$-Stirling numbers of the first kind $S_{1,r}(n+r,k+r)$ are defined by
\begin{equation}\label{eq31}
\left(x+r\right)_n=\sum_{k=0}^{n}S_{1,r}(n+r,k+r)x^k.
\end{equation}
Then
\begin{equation}\label{eq32}
\left(1+t\right)^{x+r}=\left(1+t\right)^r\sum_{k=0}^{\infty}\frac{1}{k!}\left(\log\big(1+t\big)\right)^k x^k.
\end{equation}
Further, we also have
\begin{equation}\begin{split}\label{eq33}
\left(1+t\right)^{x+r}&=\sum_{n=0}^{\infty}\left(x+r\right)_n\frac{t^n}{n!}\\
&=\sum_{n=0}^{\infty}\sum_{k=0}^{n}S_{1,r}\big(n+r,k+r\big)x^k \frac{t^n}{n!}\\
&=\sum_{k=0}^{\infty}\sum_{n=k}^{\infty}S_{1,r}\big(n+r,k+r\big)\frac{t^n}{n!}x^k.
\end{split}\end{equation}
Combining \eqref{eq32} with \eqref{eq33}, we obtain the generating function of  $S_{1,r}(n+r,k+r)$:
\begin{equation}\label{eq34}
\left(1+t\right)^r \frac{1}{k!}\left(\log\big(1+t\big)\right)^k
=\sum_{n=k}^{\infty}S_{1,r}\left(n+r,k+r\right)\frac{t^n}{n!}.
\end{equation}
\indent The central factorial numbers of the first kind $t\left(n,k\right)$ are defined by
\begin{equation}\label{eq35}
x^{[n]}=\sum_{k=0}^{n}t\left(n,k\right)x^k, \,\,\,\left(n \geq 0\right).
\end{equation}
Using \eqref{eq26} and \eqref{eq35}, we have
\begin{equation}\begin{split}\label{eq36}
\left( \frac{t}{2}+\sqrt{1+\frac{t^2}{4}}\right)^{2x}&=\sum_{n=0}^\infty x^{[n]}\frac{t^n}{n!}\\
&=\sum_{n=0}^{\infty}\sum_{k=0}^{n}t\left(n,k\right)x^k\frac{t^n}{n!}\\
&=\sum_{k=0}^{\infty}\sum_{n=k}^{\infty}t\left(n,k\right)\frac{t^n}{n!}x^k.
\end{split}\end{equation}
On the other hand, we also have
\begin{equation}\begin{split}\label{eq37}
\left( \frac{t}{2}+\sqrt{1+\frac{t^2}{4}}\right)^{2x}&=e^{2x\log \left(\frac{t}{2}+\sqrt{1+\frac{t^2}{4}}\right)}\\
&=\sum_{k=0}^{\infty}\frac{1}{k!}\left(2\log\Big(\frac{t}{2}+\sqrt{1+\frac{t^2}{4}} \Big)\right)^kx^k
\end{split}\end{equation}
By combining \eqref{eq36} with \eqref{eq37}, we get the generating function of $t\left(n,k\right)$:
\begin{equation}\label{eq38}
\frac{1}{k!}\left(2\log\Big(\frac{t}{2}+\sqrt{1+\frac{t^2}{4}} \Big)\right)^k=
\sum_{n=k}^{\infty}t\left(n,k\right)\frac{t^n}{n!}.
\end{equation}
\indent Let us define the extended $r$-central factorial numbers of the first kind as
\begin{equation}\label{eq39}
\left(x+r\right)^{[n]}=\sum_{k=0}^{n}t_{r}\left(n+r,k+r\right)x^k.
\end{equation}
Then we want to derive the generating function of  the extended $r$-central factorial numbers of the first kind.
\begin{equation}\begin{split}\label{eq40}
\left( \frac{t}{2}+\sqrt{1+\frac{t^2}{4}}\right)^{2\left(x+r\right)}&=\sum_{n=0}^\infty \left(x+r\right)^{[n]}\frac{t^n}{n!}\\
&=\sum_{n=0}^{\infty}\sum_{k=0}^{n}t_r\left(n+r,k+r\right)\frac{t^n}{n!}\\
&=\sum_{k=0}^{\infty}\sum_{n=k}^{\infty}t_r\left(n+r,k+r\right)\frac{t^n}{n!}x^k.
\end{split}\end{equation}
In addition, we also have
\begin{equation}\begin{split}\label{eq41}
\left( \frac{t}{2}+\sqrt{1+\frac{t^2}{4}}\right)^{2\left(x+r\right)}&=
\left( \frac{t}{2}+\sqrt{1+\frac{t^2}{4}}\right)^{2r}e^{2x\log \left(\frac{t}{2}+\sqrt{1+\frac{t^2}{4}}\right)}\\
&=\sum_{k=0}^{\infty}\frac{1}{k!}\left(\frac{t}{2}+\sqrt{1+\frac{t^2}{4}}\right)^{2r}\left(2\log\big(\frac{t}{2}+\sqrt{1+\frac{t^2}{4}}\big)\right)^kx^k.
\end{split}\end{equation}
Now, from \eqref{eq40} and \eqref{eq41}, we have the generating function for $t_r\left(n+r,k+r\right)$:
\begin{equation}\label{eq42}
\frac{1}{k!}\left(\frac{t}{2}+\sqrt{1+\frac{t^2}{4}}\right)^{2r}\left(2\log\Big(\frac{t}{2}+\sqrt{1+\frac{t^2}{4}}\Big)\right)^k
=\sum_{n=k}^{\infty}t_r\left(n+r,k+r\right)\frac{t^n}{n!}.
\end{equation}
\indent Finally, we want to show a recurrence relation for  the extended $r$-central factorial numbers of the first kind.
\begin{equation}\begin{split}\label{eq43}
&\sum_{k=0}^{n+1}t_{r}\left(n+1+r,k+r\right)x^k=\left(x+r\right)^{[n+1]}=\left(x+r\right)^{[n-1]}\left(\big(x+r\big)^2-\Big(\frac{n-1}{2}\Big)^2\right)\\
&=\sum_{k=0}^{n-1}t_{r}\left(n-1+r,k+r\right)x^k\left(x^2+2xr+r^2-\Big(\frac{n-1}{2}\Big)^2\right)\\
&=\sum_{k=0}^{n-1}t_{r}\left(n-1+r,k+r\right)x^{k+2}+2r\sum_{k=0}^{n-1}t_{r}\left(n-1+r,k+r\right)x^{k+1}\\
&\quad\quad\quad+\left(r^2-\Big(\frac{n-1}{2}\Big)^2\right)\sum_{k=0}^{n-1}t_{r}\left(n-1+r,k+r\right)x^k \\
&=\sum_{k=2}^{n+1}t_{r}\left(n-1+r,k-2+r\right)x^{k}+2r\sum_{k=1}^{n}t_{r}\left(n-1+r,k-1+r\right)x^{k}\\
&\quad\quad\quad+\left(r^2-\Big(\frac{n-1}{2}\Big)^2\right)\sum_{k=0}^{n-1}t_{r}\left(n-1+r,k+r\right)x^k \\
&=\sum_{k=0}^{n+1}\left\{t_{r}\left(n-1+r,k-2+r\right)+2r t_{r}\left(n-1+r,k-1+r\right)+\left(r^2-\Big(\frac{n-1}{2}\Big)^2\right)t_{r}\left(n-1+r,k+r\right)\right\}x^k.
\end{split}\end{equation}
\indent This verifies the following theorem.
\begin{theorem}
For any integers $n,k$ with $n-1 \geq k \geq 0$, we have the following recurrence relation:
\begin{equation}\begin{split}
&t_{r}\left(n+1+r,k+r\right)\\
&=t_{r}\left(n-1+r,k-2+r\right)+2r t_{r}\left(n-1+r,k-1+r\right)+\left(r^2-\Big(\frac{n-1}{2}\Big)^2\right)t_{r}\left(n-1+r,k+r\right).
\end{split}\end{equation}
\end{theorem}
\section{Conclusions and Discussions}

In recent years, quite a number of old and new special numbers and polynomials have attracted many researchers and been studied by means of generating functions, combinatorial methods, umbral calculus, differential equations, $p$-adic integrals, $p$-adic $q$-integrals, special functions, complex analysis and so on.

 In this paper, we introduced the extended $r$-central factorial numbers of the second and first kinds and the extended $r$-central Bell polynomials, and studied various properties and identities related to these numbers and polynomials and their connections. This study was done by making use of generating function techniques.

 The extended $r$-central factorial numbers of the second kind are an extended version of the central factorial numbers of the second kind and also a ‘central analogue’ of the $r$-Stirling numbers of the second; the extended $r$-central Bell polynomials are an extended version of the central Bell polynomials and also a central analogue of $r$-Bell polynomials; the extended $r$-central factorial numbers of the first kind are an extended version of the central factorial numbers of the first kind and a central analogue of the (unsigned) $r$-Stirling numbers of the first kind. All of these numbers and polynomials were studied before (see \cite{ref4,ref5,ref8,ref17}).

 As one of our next project, we would like to find some interesting applications of the numbers and polynomials introduced in this paper.

\end{document}